# Comment un décideur peut-il prendre en compte l'évolution de ses préférences au cours du temps ?

L.N. KISS, C. FONTEIX, M. CAMARGO, L. MOREL, J. RENAUD

**Résumé**

Le classement des alternatives possibles lors de la conception ou le fonctionnement d'un procédé industriel conduit à la définition d'une solution représentant le meilleurs compromis entre plusieurs objectifs contradictoires. Le choix de la solution est une décision pouvant être modélisée à l'aide d'une connaissance des préférences correspondantes. Le modèle de décision et de préférence utilisé est celui des bilans de flux. La production prise comme exemple est celle de fromages de Munster. Il s'agit de fromages très typés qui ne conviennent pas un des consommateurs non avertis. Si le décideur souhaite élargir son panel de consommateurs il doit modifier le réglage de sa ligne de production. Lors d'un changement des préférences lié à une influence extérieure la décision doit s'adapter, mais pas trop brutalement pour éviter aussi bien des oscillations de la ligne de production que pour indisposer les consommateurs. Afin de résoudre ce problème nous proposons une technique nouvelle : les bilans de flux dynamiques. Il est alors possible de prévoir l'évolution de la production. Des perspectives sont aussi proposées comme la création d'une commande prédictive multicritère.

## 1 Introduction

L'évolution du contexte socio-économique a des répercussions sur les décisions stratégiques industrielles. Ce contexte influe :

- Soit sur la valeur attribuée aux critères de décision, ce qui nécessite l'élaboration d'un modèle dynamique de calcul de ces critères,
- Soit sur la valeur des paramètres de décision, c'est-à-dire sur les préférences du décideur, ce qui nécessite l'élaboration d'un modèle d'évolution de ces paramètres,
- Soit sur l'évolution simultanée de la valeur des critères de décision et des préférences du décideur, ce qui nécessite l'élaboration de deux modèles d'évolution.

L'expression des préférences fait l'objet d'études théoriques, en particulier la décision « algorithmique ». Ces recherches se trouvent à la rencontre de deux disciplines : d'une part la théorie de la décision (KAST 2002) qui se tourne vers la modélisation et l'agrégation des préférences, et d'autre part l'algorithmique qui s'intéresse au calcul d'aide à la décision dans le risque et l'incertain (BOURDACHE 2021) (HONEL et al. 2022). Lorsque le décideur manque d'informations concordantes pour résoudre son problème, il lui arrive d'envisager une réponse sans calcul, et plutôt dynamique. L'effet de l'intuition sur la prise de décision a été étudié par (DANE et al. 2007). Mais, les informations subjectives ne sont jamais stables, et leurs fluctuations peuvent perturber la prise de décision. Alors, il arrive que le décideur cherche à se focaliser sur les données mesurables et stables. (PFIFFELMANN et al. 2013) dit : « *Notre objectif est de montrer que les différentes actions envisageables ne sont pas analysées objectivement mais subjectivement par le décideur. En effet, grâce à l'économie expérimentale, il est possible de mettre en évidence un certain nombre de biais de comportement chez le décideur tels que l'optimisme, le conservatisme ou l'aversion aux pertes.* » De même, (FARQUHAR et al. 1989) indiquent : « *The concept of preference intensity has been criticized over the past sixty years for having no substantive meaning. Much of the controversy stems from the inadequacy of*

*measurement procedures.* » En fait, l'évolution du décideur a toujours été guidée par les différentes transformations sociétales, technologiques, et environnementales apparues au cours du temps. Ces évolutions ont contraint les entreprises à modifier en permanence leurs offres et à innover. (LA BRUSLERIE et al. 2012) précisent : « *La préférence intemporelle est une hypothèse de comportement. Elle se réfère à une préférence pour l'utilité immédiate par rapport à une utilité future.* » et « *Le temps est à la fois une ressource économique, une dimension du choix et un lieu d'expression de l'utilité.* » La plupart des décisions sont complexes et incertaines, (BOUQUET 2013) dit : « *Aussi beaucoup de réflexions portent sur la rationalité et l'irrationalité de la décision, la décision complexe, les paradoxes de la décision, et la temporalité de la décision y est régulièrement évoquée.* »

Dans ce papier la section 2 présentera l'application industrielle et la section 3 la nouveauté de ce travail : les bilans de flux dynamiques. La section 4 montrera comment adapter la production à de nouveaux consommateurs et la section 5 appliquera les bilans de flux dynamiques à un échantillon de 5 produits. La section 6 détaillera des perspectives et la section 7 nos conclusions.

## 2 Application industrielle

Afin de tester notre approche un ensemble de données a été recueilli au sein d'une entreprise fromagère (BOUCHOUX 2002). Il s'agissait d'une fabrication de munsters de 17 cm de diamètre, et 47 fromages ont été goûtés par deux à six professionnels de l'analyse sensorielle sur une période de trois mois. Dix-huit critères ont été notés de 0 à 7 (en nombre décimal) pour chaque fromage et l'écart à une cible a été calculée pour chaque critère. Une moyenne a été déduite pour chaque critère, puis certains critères ont été agrégés pour qu'il n'en reste que quatre. Par exemple, le « grand » critère correspondant à l'aspect extérieur regroupe 5 sous critères. Il en est de même de la texture et du goût, ou arôme, par contre l'odeur, ou le parfum, n'en regroupe que trois. Ainsi nous obtenons 4 critères par fromage : l'aspect, noté C1, l'odeur, notée C2, la texture, notée C3 et le goût noté C4. L'écart à la cible est moyenné et transformé à l'aide d'une quantification floue présentée dans (RENAUD et al. 2008). La valeur prise en compte est donc comprise entre 0 et 1. Un échantillon de 5 fromages a été retenu pour la présente étude, il s'agit de l'échantillon E1 de (RENAUD et al. 2008). Le tableau 1 montre les notes retenues pour chaque fromage et chaque critère. Le numéro de chaque fromage est un code donné par l'entreprise.

| *Ech.1* | C1 | C2 | C3 | C4 |
|---|---|---|---|---|
| **613** | 0,62093 | 0,70547 | 0,734 | 0,99189 |
| **2573** | 0,8907 | 0,85185 | 0,666 | 0,54054 |
| **292** | 0,81395 | 0,97002 | 0,4 | 0,33784 |
| **162** | 0,77442 | 0,82363 | 0,734 | 0 |
| **3062** | 0,5814 | 0,17637 | 0,7 | 0,67568 |

Tableau 1 : Valeur moyenne de chaque critère pour chaque fromage

Ce code permet à l'entreprise de connaitre précisément les conditions de production du fromage correspondant. Notons que les fromages notés de 613 à 3062 sont présentés dans l'ordre du choix du décideur, du meilleur au moins bon. En effet, la production est complexe et de nombreux facteurs sont difficilement maitrisables, comme par exemple la qualité et la composition du lait utilisé. Ces éléments sont liés à la saison, à l'alimentation des animaux et à l'exploitant du domaine laitier.

## 3 Le modèle de décision multicritère : les bilans de flux dynamiques

Une prise de décision peut être un processus complexe (BERARD 2009), qui doit être performant en situation d'incertitude (OUZOUNOVA 2005). Dans le répertoire des modèles de décision multicritère et de préférences humaines (ENJOLRAS et al. 2022), nous avons opté pour le concept des bilans de flux (KISS et al. 1994) (DEROT et al. 1997). Ce modèle a été appliqué dans (KISS et al. 2002) (FONTEIX et al. 2004) (RENAUD et al. 2007). Dans ce cas le décideur est confronté à un problème de décision multicritère, et les bilans de flux doivent classer les alternatives envisagées de la préférée à la moins intéressante. Pour cela le décideur indique ses préférences en fonction de sa connaissance du procédé, puis l'algorithme compare les alternatives deux à deux. Pour effectuer cette comparaison le décideur doit définir ses préférences à travers plusieurs seuils. Avant cela, il doit préciser les poids de chaque critère, dépendant de leur importance relative sous une forme normalisée, $w_k$ pour le $k^{ième}$ critère sur $n$ tels que :

$$\sum_{k=1}^{k=n} w_k = 1 \qquad (1)$$

De plus, le décideur doit indiquer la valeur $q_k$ du seuil d'indifférence attribué au critère $k$. Ce seuil veut dire qu'en dessous le décideur n'est pas capable de distinguer la valeur du critère $k$ d'une alternative de celle d'une autre. C'est à dire que si la différence entre le critère $k$ de deux alternatives est inférieure ou égale au seuil d'indifférence (en valeur absolue) alors les deux alternatives sont aussi intéressantes l'une que l'autre. Par contre, si cette différence est supérieure ou égale au seuil de préférence $p_k$, l'alternative ayant la plus forte valeur du critère $k$ est préférée à l'autre. Dans le cas où cette différence est supérieure au seuil de véto $v_k$, l'alternative de moins bonne (plus faible) valeur du critère $k$ est considérée comme trop mauvaise et est pénalisée quel que soit la valeur des autres critères. Pour chaque critère, ces trois seuils doivent respecter la contrainte suivante :

$$0 \leq q_k \leq p_k \leq v_k \qquad (2)$$

Tous les critères étant à maximiser, posons $\Delta_k(i,j) = C_k(i) - C_k(j)$ la différence indiquée plus haut entre l'alternative $i$ et l'alternative $j$ pour le critère $C_k$. Cette différence va donner deux valeurs, une de concordance (entre les critères) et une de discordance (si une alternative est pénalisée comme indiqué plus haut) comme le montre la figure 1.

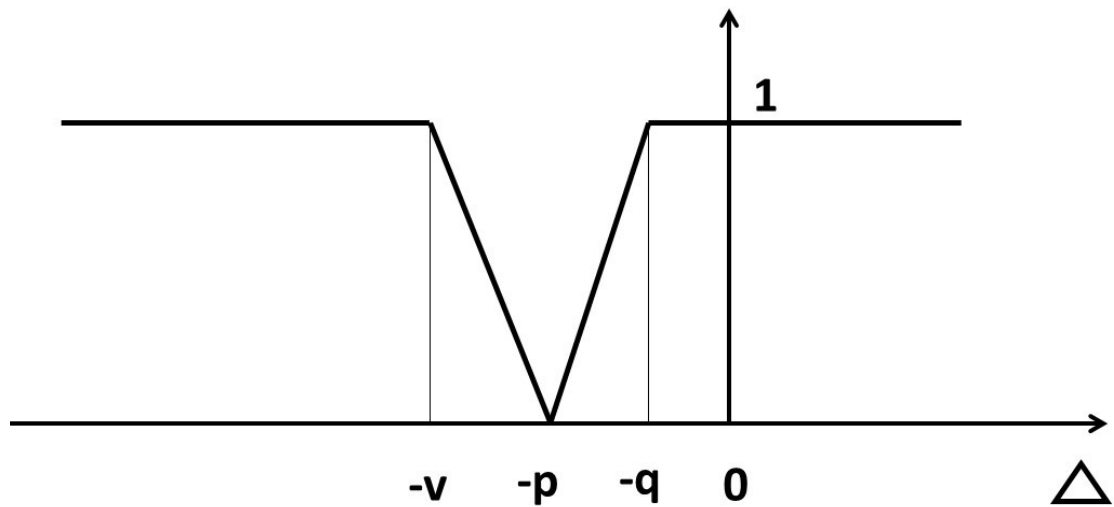

Figure 1 : courbes de concordance et de discordance

La préférence d'une alternative sur une autre est dite concordante si nous nous trouvons sur le plateau de droite égal à 1. Entre les seuils $-p_k$ et $-q_k$ la concordance est partielle et suit l'expression :

$$c_k(i,j) = 0 \quad si \quad i = j$$

$$sinon \quad c_k(i,j) = \begin{cases} 1 & si & \Delta_k(i,j) \geq -q_k \\ \frac{\Delta_k(i,j)+p_k}{p_k-q_k} & si & -p_k < \Delta_k(i,j) < -q_k \\ 0 & si & \Delta_k(i,j) \leq -p_k \end{cases} \quad (3)$$

La discordance correspond au plateau de gauche égal à 1, figure 1, et son expression est :

$$D_k(i,j) = \begin{cases} 1 & si & \Delta_k(i,j) \leq -v_k \\ \frac{-\Delta_k(i,j)-p_k}{v_k-p_k} & si & -v_k < \Delta_k(i,j) < -p_k \\ 0 & si & \Delta_k(i,j) \geq -p_k \end{cases} \quad (4)$$

En utilisant la concordance et la discordance il est possible de définir un degré de surclassement d'une alternative sur une autre grâce à (FONTEIX et al. 2004) :

$$\sigma(i,j) = \left(\sum_{k=1}^{n} w_k c_k(i,j)\right)\left(\prod_{k=1}^{n}(1 - D_k(i,j)^3)\right) \quad (5)$$

Finalement, le flux net de supériorité de l'alternative $i$ sur l'ensemble des $m$ alternatives est :

$$\varphi^+(i) = \sum_{j=1}^{m} \sigma(i,j) \quad (6)$$

De même, nous définissons le flux net d'infériorité de l'alternative $i$ sur les autres par :

$$\varphi^-(i) = \sum_{j=1}^{m} \sigma(j,i) \quad (7)$$

Le score usuellement utilisé pour classer les alternatives de la première à la dernière est :

$$s(i) = \varphi^+(i) - \varphi^-(i) \quad (8)$$

Cependant, nous sommes ici dans le cas où rien ne change au cours du temps, si par contre une évolution s'avère nécessaire parce que les paramètres, seuils et poids, doivent changer pour s'adapter à l'environnement de la décision l'équation (8) s'avère inadéquate. En fait, le changement peut être subit par le décideur s'il s'agit d'une modification des critères $C_k(i)$. Il peut aussi être anticipé par ce dernier s'il est amené à choisir de nouveau paramètres, poids et deuils, pour s'adapter à une tendance constatée chez ses clients. Ainsi, l'évolution du contexte socio-économique a des répercussions sur les décisions stratégiques industrielles. Comme nous l'avons dit en introduction, le contexte influe soit sur la valeur attribuée aux critères de décision, soit sur la valeur des paramètres de décision, soit sur l'évolution simultanée de la valeur des critères de décision et des préférences du décideur.

Dans tous les cas, il faut envisager un nouveau modèle d'aide à la décision, les bilans de flux dynamiques, car l'accompagnement du changement ne doit pas être brutal, il doit même amortir les fluctuations possibles de l'environnement de la décision. Si nous prenons l'exemple des valeurs boursières, celles-ci fluctuent selon deux modes distincts : des oscillations à haute fréquence qui se superposent à des oscillations à basse fréquences. Ces dernières représentent l'évolution effective du marché alors que les premières ne correspondent qu'à sa nervosité. Il en est de même pour le contexte socio-économique des décisions industrielles, or un décideur doit prendre en compte les fluctuations à basse fréquence, qui représentent les tendances lourdes de son environnement, et ignorer autant que possible les fluctuations à haute fréquence, qui ne sont que des épiphénomènes non durables. Le filtre passe bas est un outil largement utilisé en économie pour décomposer une série chronologique dans sa tendance et ses composants cycliques. En séparant la tendance à long terme des fluctuations à court terme, les décideurs peuvent obtenir des informations précieuses sur les conditions économiques sous-jacentes et prendre des décisions éclairées. Une solution simple à ce problème est d'utiliser un filtre numérique du 1$^{er}$ ordre dans le calcul des scores par alternative. Ainsi nous remplacerons l'équation statique (8) par la forme dynamique :

$$s(i) + \tau \frac{ds(i)}{dt} = \varphi^+(i) - \varphi^-(i) \qquad (9)$$

Où $t$ est le temps et $\tau$ un temps d'amortissement choisi pour répartir l'évolution rapide de $\varphi^+(i) - \varphi^-(i)$ entre $s(i)$ et $\tau \frac{ds(i)}{dt}$ de façon à réduire l'impact de la variation de $\varphi^+(i) - \varphi^-(i)$ sur $s(i)$. Le résultat est un vecteur des scores qui évolue au cours du temps et qui peut conduire à une évolution effective de la décision choisie par le décideur au fil du temps, car le classement des alternatives peut changer. Il s'agit là d'une vraie décision du décideur humain qui peut choisir une alternative légèrement déclassée si cela lui évite des changements trop fréquents et importants au sein de l'entreprise. Le choix de la valeur de $\tau$ est important pour bien réduire l'impact des fluctuations rapides du contexte socio-économique, et ce paramètre peut être considéré comme un paramètre supplémentaire de décision.

En pratique le calcul des scores ne se fait pas à chaque instant, mais tous les $\Delta t$ :

$$s(i, t + \Delta t) + \frac{\tau}{\Delta t}\big(s(i, t + \Delta t) - s(i, t)\big) = \varphi^+(i, t) - \varphi^-(i, t) \qquad (10)$$

Si nous appliquons l'intégration d'EULER, ainsi l'équation (9) correspondait à un filtre analogique, que nous avons transformé dans (10) en un filtre numérique plus facile à utiliser dans le contexte présent. Dans ce cas, nous poserons $\alpha = \frac{1}{1 + \frac{\tau}{\Delta t}}$, ce qui nous conduit à écrire le filtre numérique sous la forme :

$$s(i, t + \Delta t) = (1 - \alpha)s(i, t) + \alpha\big(\varphi^+(i, t) - \varphi^-(i, t)\big) \qquad (11)$$

Il s'agit là de la forme définitive des bilans de flux dynamiques. Il faut calculer les scores de toutes les alternatives tous les $\Delta t$, en tous cas si quelque chose a changé, soit dans le calcul des critères, soit dans le choix des préférences du décideur, c'est-à-dire dans le paramétrage des bilans de flux.

**4 Adapter la production aux souhaits des consommateurs**

La production de fromages de l'entreprise est destinée à être vendue aux consommateurs. Or l'évolution du comportement des consommateurs a toujours été guidée par les différentes transformations sociétales, technologique et environnementales apparues au fil du temps. Ces évolutions ont contraint les entreprises à modifier en permanence leurs offres et à innover (LADWEIN 2003) (BREE 2017) (ALAMI et al. 2023). Reconnaître ces évolutions est une chose, mais les analyser et les interpréter de manière efficace et optimale en est une autre (GALIBERT et al. 1986) (ROBIN et al. 2007). L'aspect environnemental, social et sociétal représente aujourd'hui une motivation importante des consommateurs (TREMBLAY 1994) (ANDALOUSSI 2002) (KALLEL 2007) (BASTI 2010) (CHOW-YING 2013) (CLAUZEL et al. 2016) (WOOD 2018). Le consommateur réfléchit, devient prudent et a tendance à se reporter sur les marques qui lui semblent les plus fiables et qui offrent le meilleurs rapport qualité/prix. D'ailleurs, cela assure le succès des produits qui respectent l'authenticité du terroir. C'est sur ce point précis que nous allons nous appuyer pour envisager le changement d'avis du décideur.

Les notes de 0 à 7 mises par les goûteurs se sont avérées précises et reproductibles à moins de 0,1 près. Les objectifs par critère sont stables et indépendants du choix du décideur quant au fromage produit. Ce fait se transmet aux moyennes retenues après comparaison des résultats donnés par chaque goûteur. En effet, les goûteurs travaillent séparément en simultané dans des box isolés phoniquement. Ainsi, les distances à l'objectif entre 0 et 1 obtenues selon la procédure indiquée dans (RENAUD et al. 2008) sont précises à moins de 0,01. Ceci permet de définir les seuils de façon que :

| Seuil | Valeur |
|---|---|
| Indifférence | 0 |
| Préférence | 0,1 |
| Véto | 0,3 |

Tableau 2 : Valeur donnée aux seuils

Jusqu'à présent, le décideur avait choisi de produire un fromage très typé, convenant à des consommateurs avertis. Il s'agissait donc de fromages traditionnels correspondant à une niche de consommation. Les poids des 4 critères étaient donc :

| Critère | Type | Poids |
|---|---|---|
| C1 | Aspect | 0,1 |
| C2 | Odeur | 0,4 |
| C3 | Texture | 0,1 |
| C4 | Goût | 0,4 |

Tableau 3 : Poids initiaux des critères

Afin de mieux voir le fonctionnement du calcul des notes données à chaque fromage nous réécrirons les équations de (5) à (8) :

$$\left.\begin{array}{l}\sigma(i,j) = \sum_{k=1}^{n} w_k \left(c_k(i,j) \prod_{k=1}^{n}(1 - D_k(i,j)^3)\right) = \sum_{k=1}^{n} w_k \, \omega_k(i,j) \\ \varphi^+(i) = \sum_{j=1}^{m} \sigma(i,j) = \sum_{k=1}^{n} w_k \sum_{j=1}^{m} \omega_k(i,j) = \sum_{k=1}^{n} w_k \, \mu_k^+(i) \\ \varphi^-(i) = \sum_{j=1}^{m} \sigma(j,i) = \sum_{k=1}^{n} w_k \sum_{j=1}^{m} \omega_k(j,i) = \sum_{k=1}^{n} w_k \, \mu_k^-(i)\end{array}\right\} \quad (12)$$

C'est-à-dire :

$$\left.\begin{array}{l}\mu_k^+(i) = \sum_{j=1}^{m} \omega_k(i,j) = \sum_{j=1}^{m} c_k(i,j) \prod_{k=1}^{n}(1 - D_k(i,j)^3) \\ \mu_k^-(i) = \sum_{j=1}^{m} \omega_k(j,i) = \sum_{j=1}^{m} c_k(j,i) \prod_{k=1}^{n}(1 - D_k(j,i)^3) \\ s(i) = \sum_{k=1}^{n} w_k \left(\mu_k^+(i) - \mu_k^-(i)\right)\end{array}\right\} \quad (13)$$

Le tableau des résultats du calcul des bilans de flux est :

| *Ech.1* | C1 | C2 | C3 | C4 |
|---|---|---|---|---|
| **613** | 1 | 1 | 2,76214551 | 2,76214551 |
| **2573** | 2,90774712 | 1,62434345 | 1,59214939 | 1,6157292 |
| **292** | -0,91317107 | 0,3702326 | -1,39773957 | -1,39773957 |
| **162** | -1 | -1 | -1,30013514 | -1,98013514 |
| **3062** | -1,99457605 | -1,99457605 | -1,6564202 | -1 |

Tableau 4 : Valeurs de $\mu_k^+(i) - \mu_k^-(i)$ en fonction du fromage et du critère

Avec l'application des poids et seuils précédents à l'échantillon E1, nous obtenons les notes suivantes :

| *Ech.1* | Note $s(i)$ |
|---|---|
| **613** | 1,88107276 |
| **2573** | 1,74601871 |
| **292** | -0,64209385 |
| **162** | -1,42206757 |
| **3062** | -1,56293005 |

Tableau 5 : Score de chaque fromage avec les poids du tableau 3

Nous constatons que le classement obtenu par les notes, ou scores, données par les bilans de flux correspond bien à celui du décideur.

Supposons maintenant que suite à une étude marketing le décideur souhaite augmenter sa production en se conformant plus au goût des consommateurs étrangers, ou de ceux n'aimant pas les fromages trop typés. Alors il faudrait produire des fromages selon les critères :

| Critère | Type | Poids |
|---|---|---|
| C1 | Aspect | 0,4 |
| C2 | Odeur | 0,1 |
| C3 | Texture | 0,4 |
| C4 | Goût | 0,1 |

Tableau 6 : Poids des critères pour un fromage moins typé

Les seuils n'étant pas modifiés, nous pouvons toujours utiliser les résultats données tableau 4, donc à l'aide des tableaux 4 et 6 nous pouvons calculer les nouveaux scores des fromages de E1 :

| Ech.1 | Note $s(i)$ |
|---|---|
| **613** | 1,88107276 |
| **2573** | 2,12396587 |
| **292** | -1,02711495 |
| **162** | -1,21806757 |
| **3062** | -1,75985611 |

Tableau 7 : Score de chaque fromage avec les poids du tableau 6

Nous constatons une inversion du classement des 2 premiers fromages. Cela veut dire que pour produire des fromages moins typés il convient d'appliquer les conditions de production ayant conduit à l'obtention du fromage 2573.

**5 Application à l'échantillon E1**

Notre objectif est de passer d'une production de fromages typés à une production de fromages moins typés. Le changement brutal des conditions opératoires de production peut engendrer des difficultés, voire des oscillations. D'autre part le changement brutal des qualités fromagères du produit peut indisposer les consommateurs. De ce fait un changement progressif est généralement préférable. Pour cela nous appliquerons l'équation (9) avec trois valeurs différentes de $\alpha = \frac{1}{1+\frac{\tau}{\Delta t}}$, c'est-à-dire 0,1 puis 0,3 et 0,5. Plus $\alpha$ est grand plus vite l'objectif sera atteint. La figure 2 montre l'évolution des scores des fromages 613 et 2573 pour les trois valeurs de $\alpha$ indiquées plus haut :

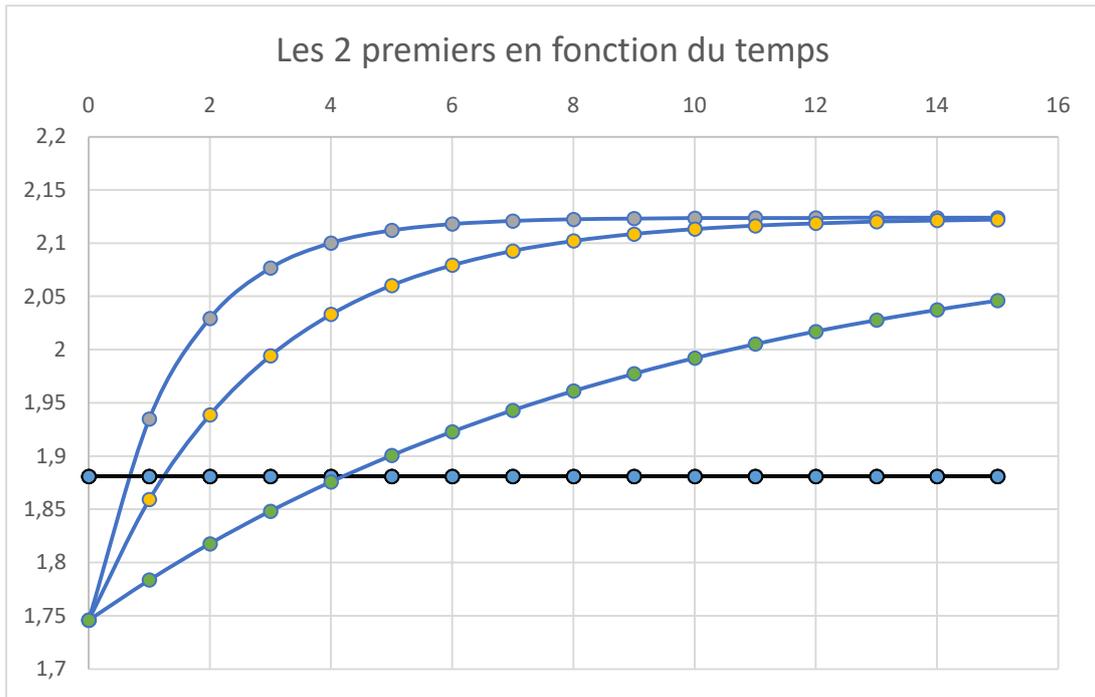

Figure 2 : Evolution des scores des fromages 613 et 2573 au cours du temps

Nous constatons que la note attribuée au fromage 613 n'est pas modifiée par un changement de valeur des poids. Rappelons que c'était le fromage préféré du décideur, cela peut donc s'expliquer par une excellente connaissance de sa ligne de production. Il savait que le produit 613 était très équilibré. Par contre les trois courbes correspondant au fromage 2573 coupent celle du produit 613, la plus haute correspondant à $\alpha = 0{,}5$ et la plus basse à $\alpha = 0{,}1$. Les courbes se croisent, il y a donc inversion dans le classement, en tant que fromage moins typé, 2573 est le nouveau produit préféré. En ce qui concerne les trois derniers fromages du classement, nous ne constatons aucun changement de leur classement.

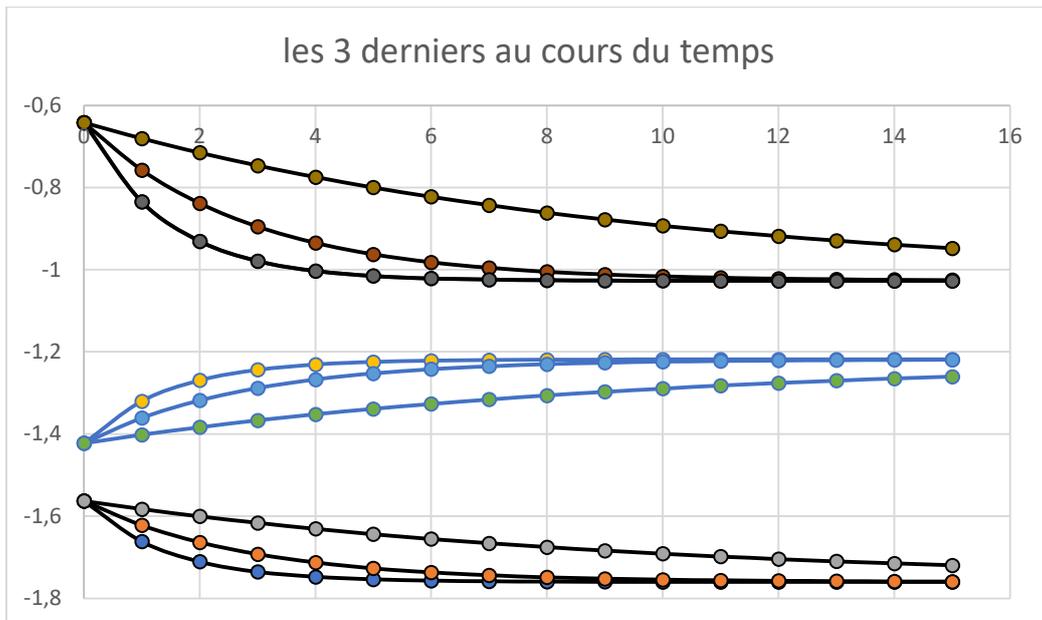

Figure 3 : Evolution des scores des fromages 292, 162 et 3062 au cours du temps

Les figures 2 et 3 montrent comment le classement de l'échantillon E1 change au cours du temps. L'unité de temps est $\Delta t$, le pas de calcul du filtre que nous ne préciserons pas ici. La taille de l'échantillon est faible, 5, ce qui rend la visualisation de l'évolution du classement par les bilans de flux dynamiques facile à comprendre, mais les choses seraient moins commodes si le nombre de fromages était beaucoup plus grand. Cependant, le calcul ne serait pas plus compliqué, seule la visualisation serait moins simple.

**6 Perspectives**

Le classement des alternatives possibles (choix des procédés ou matériels, répartition de la demande sur plusieurs installations…) lors de la conception ou du fonctionnement d'un ou plusieurs procédés des industries chimiques ou parachimiques (alimentaires, pharmaceutiques, cosmétiques…) conduit à la définition d'une solution représentant le meilleur compromis entre plusieurs objectifs contradictoires, comme la rentabilité, la qualité des produits ou leurs propriétés d'usage, leur coût écologique, etc… Ces critères, conduisant généralement à des solutions à critères contradictoires, nécessitent la détermination préalable d'un espace de Pareto, tel que tout point ne se trouvant pas dans ce domaine soit une solution totalement prohibée, puisqu'il est certain qu'il en existe de bien meilleures (VIENNET et al. 1995) (VIENNET et al. 1996). Il est donc impératif que la solution retenue reste dans ce domaine de Pareto, voire une partie de ce dernier considérée comme la zone d'intérêt industriel, sous l'effet des perturbations naturelles, comme la variabilité des matières premières ou du procédé lui-même. Ainsi, le classement effectué doit intégrer la robustesse technique, c'est-à-dire doit peu évoluer sous l'influence de changements, d'incertitudes ou de perturbations. A partir de la variance ou de l'encadrement des perturbations, une distance entre une décision candidate et le bord de la zone d'intérêt est définie. Cette dernière est ensuite introduite comme critère à maximiser supplémentaire ajouté aux objectifs déjà existants.

Il est apparu que les outils du Génie Décisionnel, ou de l'Analyse Multicritère, étaient très complémentaires à l'optimisation multicritère au sens de Pareto. Il s'agit d'avancer en deux phases, d'abord réaliser cette optimisation multicritère qui fournit le domaine de Pareto, c'est-à-dire les compromis, et ensuite classer les alternatives à l'aide du Génie Décisionnel (MASSEBEUF et al. 1999) (KISS et al. 2002). La notion de robustesse décisionnelle se limite généralement à une étude de sensibilité de la décision aux paramètres du modèle de préférence. En pratique, dans le secteur industriel, il est préférable de rechercher la décision de meilleur rang, la moins sensible à une variation des paramètres de décision, c'est à dire la solution qui minimise le rang le plus élevé, obtenu dans l'espace de variation des paramètres. Cette robustesse décisionnelle doit permettre à la solution envisagée de rester une bonne décision, même si l'évolution du contexte politique ou économique oblige le décideur à revoir un peu son échelle de préférence. En fait, trois types de robustesse ont été définis pour chaque alternative acceptable, la robustesse technique vue plus haut, liée aux aléas de la production, la robustesse décisionnelle, liée à une évolution possible du contexte de la décision (technico-économique, politique…), et la robustesse procédurale, liée à l'influence de la technique de classement utilisée. Cette dernière présente un intérêt certain, car si nous utilisons plusieurs modèles de décision nous obtenons plusieurs classements des compromis (alternatives Pareto efficaces). Alors il est possible de proposer au choix du décideur plusieurs solutions différentes et bien classées (autant que de modèles de préférence). Mieux, il est possible de comparer les différents classements et de les considérer comme des nouveaux critères à minimiser. Alors en dégageant le domaine de Pareto des classements on réduit le nombre d'alternatives à proposer au décideur et on est sûr qu'il s'agit des meilleures (MUNIGLIA et al. 2001).

Nous avons montré dans les figures 2 et 3 que les bilans de flux dynamiques étaient possibles. La prise en compte d'objectifs multiples dans la définition des actions à réaliser en temps réel, sur un procédé, constitue une nouveauté scientifique, mais aussi une suite logique à l'introduction de la dynamique dans l'optimisation multicritère. Il s'agit là de développer une **commande prédictive multicritère** par fusion de notions issues de l'Automatique, du Génie des Procédés et du Génie Décisionnel. Dans sa thèse (CLIVILLE 2004) dit « *Si toutes les relations cause-effet peuvent être identifiées et quantifiées, le système piloté est totalement contrôlable et le pilotage est analogue à la commande de l'automatique. Dans le cas général, le système piloté est trop complexe, et seule une partie de ces relations est formellement identifiée. Dans ce cas, le système de pilotage génère un ou plusieurs plans d'action pour atteindre les objectifs. Il faut alors prendre des décisions (choisir, trier, classer parmi un ensemble de solutions).* » Le contrôle de la dynamique de la décision (décision variable dans le temps, de façon parfois importante, tout comme une commande) doit permettre un filtrage « numérique » de cette dernière, et une bonne robustesse (rejet de perturbations) du système. La commande multicritère utilise un modèle du procédé et un outil de génie décisionnel. Cet outil est un modèle de décision humaine qui doit être adapté au problème posé. Avant toute décision multicritère, le décideur exprime ses préférences qu'il faut évaluer ou mesurer. Les préférences ne sont pas toujours faciles à exprimer, et leur modélisation a parfois une structure très complexe. Deux éléments « clés » dans la problématique multicritère sont le choix des critères et leurs importances ou leurs « poids ». Evaluer les poids des critères reste une opération délicate. Nous proposons deux nouvelles approches :

- par classement d'un échantillon de cas par le décideur. Les scores calculés par la méthode d'agrégation pour ces cas devront respecter ce classement. Moyennant une hypothèse d'équipartition de ces scores, il est possible d'identifier les paramètres du modèle de décision,
- par la notation par le décideur de chaque cas d'un échantillon. La minimisation de l'écart entre le score calculé par la méthode d'agrégation pour chaque cas et celui fournit par le décideur permet d'identifier les paramètres de décision (poids, seuils...). Cependant une problématique scientifique importante reste inexplorée : le choix de l'échantillon optimal à présenter aux décideurs. La solution envisagée consiste à adapter à ce problème la technique d'optimisation des plans d'expériences (D-optimalité).

La commande prédictive non-linéaire, ou commande optimale à horizon fuyant, permet à un système de suivre une trajectoire prédéfinie dans l'espace d'état (consigne), ou à optimiser un critère donné. La prise en compte d'objectifs multiples dans la définition des actions à réaliser en temps réel, sur un procédé, constitue une nouveauté scientifique, mais aussi une suite logique à l'introduction de la dynamique dans l'optimisation multicritère. Il s'agit là de développer une **commande prédictive multicritère** par fusion de notions issues de l'Automatique, du Génie des Procédés et du Génie Décisionnel.

**7 Conclusions**

Les résultats des deux figures précédentes montrent qu'il est possible d'appliquer la méthode des bilans de flux dans le cas d'un changement des poids par critère. Nous pensons donc pouvoir affirmer que ceci est généralisable aux changements de définition des critères et de valeur des seuils. Ainsi, nous disposons d'un nouvel outil multicritère, les bilans de flux dynamiques, équation (11). A l'aide de cet outil il est possible de changer progressivement de décision en fonction de l'évolution du système et de son contexte. Dans les perspectives nous avons souligné qu'il fallait envisager trois types différents de robustesse : les robustesses technique, décisionnelle et procédurale. Chacune d'entre elle a son intérêt.

Par ailleurs, les bilans de flux dynamiques nous conduisent naturellement à envisager une commande multicritère, voire une commande prédictive multicritère. Tout un champ de recherche s'ouvre dans cette direction. De plus, l'identification paramétrique, ce que les informaticiens appellent calage de code, doit pouvoir s'appliquer aux modèles de préférence multicritères. Ceci nous fait penser à un moyen de sélectionner un ensemble de données parmi les compromis afin d'obtenir les valeurs des paramètres du modèle (poids et seuils pour les bilans de flux) à partir de mesures. Cela veut dire que nous devons demander au décideur de classer ou noter (à sa convenance) l'ensemble des données issues du domaine de Pareto. Afin de minimiser l'erreur sur la valeur des paramètres il convient d'utiliser une technique issue des plans d'expériences (D-optimalité par exemple) pour extraire l'ensemble d'identification de tous les compromis utilisables.